\documentclass[preprint,12pt]{elsarticle}

\usepackage{amsmath}
\usepackage{amssymb}
\usepackage{amsthm}
\usepackage{epstopdf}
\usepackage{subfigure}
\usepackage{float}
\usepackage{tikz}

\newtheorem{theorem}{Theorem}[section]
\newtheorem{lemma}[theorem]{Lemma}

\newtheorem{corollary}[theorem]{Corollary}
\journal{Journal of Combinatorial Theory Series B}

\begin{document}

\begin{frontmatter}



\title{Nonpositive Eigenvalues of the Adjacency Matrix and Lower Bounds for Laplacian Eigenvalues\tnoteref{label1}}
\tnotetext[label1]{This research supported by National Science Foundation grant DMS-0751964. This work was also performed in part using
computational facilities at the College of William and Mary which were provided with the assistance of the National Science Foundation, the Virginia Port Authority, Sun Microsystems, and Virginia's Commonwealth Technology Resource Fund.}
\author[ZBC]{Zachary B. Charles}
\author[MF]{Miriam Farber}
\author[CRJ]{Charles R. Johnson}
\author[LKS]{Lee Kennedy-Shaffer}
\address[ZBC]{University of Pennsylvania, Department of Mathematics, Philadelphia, PA 19104, USA. zcharles@sas.upenn.edu}
\address[MF]{Technion -- Israel Institute of Technology, Department of Mathematics, Haifa, IL-32000, Israel, tel 972-542345759, fax 972-778840189.\\ miriamf@tx.technion.ac.il (corresponding author)}
\address[CRJ]{College of William and Mary, Department of Mathematics, Williamsburg, VA 23187-8795, USA. crjohnso@math.wm.edu}
\address[LKS]{Yale University, Department of Mathematics, New Haven, CT, 06520, USA. lee.kennedy-shaffer@yale.edu}

\begin{abstract}
Let $NPO(k)$ be the smallest number $n$  such that the adjacency matrix of any undirected graph with $n$ vertices or more has at least $k$ nonpositive eigenvalues. We show that $NPO(k)$ is well-defined and prove that the values of $NPO(k)$ for $k=1,2,3,4,5$ are $1,3,6,10,16$ respectively. In addition, we prove that for all $k \geq 5$, $R(k,k+1) \ge NPO(k) > T_k$, in which $R(k,k+1)$ is the Ramsey number for $k$ and $k+1$, and $T_k$ is the $k^{th}$ triangular number. This implies new lower bounds for eigenvalues of Laplacian matrices: the $k$-th largest eigenvalue is bounded from below by the $NPO(k)$-th largest degree, which generalizes some prior results.
\end{abstract}

\end{frontmatter}

\pagebreak

\section{Introduction}
Given a number of vertices $n$, how many nonpositive eigenvalues \textbf{must occur in any} $n-\text{by}-n$ adjacency matrix ? Equivalently, we may ask to bound the inertia of the adjacency matrix (the number of positive, negative, and zero eigenvalues). The answer to this question is not only interesting on its own, but also relates to other algebraic graph theory questions. For example, the independence number of a graph is bounded from above by the number of nonpositive eigenvalues \cite{Cvetkovic}. Importantly to us is the close connection between the number of nonpositive eigenvalues of the adjacency matrix and bounds for Laplacian eigenvalues through the diagonal entries of the Laplacian matrix. This was our motivation, and is described herein.\\

It is well-known that majorization provides a complete description of the relationship between the possible spectra and diagonal entries of an Hermitian matrix \cite{CJ}. For certain subclasses of Hermitian matrices, however, additional inequalities that restrict this relationship may occur. We consider the Laplacian matrix of a graph, whose eigenvalues have been widely studied \cite{J. Cheeger, M. Fiedler, Merris,Ed. Y}. An example of a connection between the Laplacian eigenvalues and the degrees of the vertices in a graph is given in the following theorem \cite{Brouwer_Haemers}:
\begin{theorem}
    Let $G$ be a finite, simple, unweighted graph on $n$ vertices. If $G$ is not $K_{i} \cup (n-i)K_1$, then $\lambda_i(G) \geq d_i(G)+2-i$, in which $\lambda_i(G)$ is the $i^\text{th}$ largest Laplacian eigenvalue, and $d_i(G)$ is the $i^\text{th}$ largest degree of $G$.
\end{theorem}

We ask what other, possibly simpler, relations exist. Specifically, we are interested in the smallest integer $j$ such that for any graph on $m \ge j$ vertices, the $k^\text{th}$ largest eigenvalue of the Laplacian is at least the $j^\text{th}$ largest degree in the graph. There are no such relationships for weighted graphs, as shown in \cite{C_F_J_K}. Interestingly, there are for unweighted graphs.
Here, we first study the number of nonpositive eigenvalues of the adjacency matrix of a graph. Certain named graphs and families of graphs, such as the Kneser graphs, Paley graphs, Petersen graph and Clebsh graph, play an important role in many of our results. Then, using these results, we derive new lower bounds for Laplacian eigenvalues of a graph.

\vspace{6pt}

\section{Definitions and Key Lemmas}
Let $G$ be a simple undirected graph on $n$ vertices. Without loss of generality, we may label its vertices so that the $i$-th vertex has degree $d_i$, and $d_1 \ge d_2 \ge \cdots \ge d_n$. The smallest degree of a vertex in $G$ is denoted by $\delta(G)$, and the complement of $G$ is denoted by $\overline{G}$. The adjacency matrix of the graph $G$ is denoted by $A(G)$, and its Laplacian matrix , $L(G)$, is defined as $L(G)=diag(d_1, d_2, \dots, d_n)-A(G)$. We denote by $\lambda_1 \ge \lambda_2 \ge \cdots \ge \lambda_n = 0$ the eigenvalues of $L(G)$.\\
The Schur complement will prove useful in investigating the eigenvalues of submatrices. Consider a block matrix $M=\left[\begin{matrix} M_{11} & M_{12} \\ M_{21} & M_{22}\end{matrix}\right]$ such that $M_{22}$ is invertible. The \emph{Schur complement} (\cite{CJ}) of $M_{22}$ is given by the matrix
\begin{center}
$M/M_{22}=M_{11}-M_{12}M^{-1}_{22}M_{21} $.\\
\end{center}
The inertia of a matrix $M$ is the ordered triple $i(M)=(i_+(M),i_-(M),i_0(M))$, in which $i_+(M)$,$i_-(M)$ and $i_0(M)$ are the numbers (counting multiplicity) of positive, negative, and zero eigenvalues of $M$, respectively \cite{CJ}.
The following lemma will be crucial in our use of the Schur complement.
\begin{lemma} \label{lemm:schur} (\cite{inertia})
For an $n-\text{by}-n$ Hermitian block matrix $M$ partitioned as above, $i(M)=i(M_{22})+i(M/M_{22}) $.
\end{lemma}

Also important will be the Interlacing Theorem (\cite{CJ}) , given in our notation as follows:
\begin{theorem}\label{thm:inter} Let $A\in M_n$ be a given Hermitian matrix, and let $B \in M_{n-1}$ be a principal submatrix of $A$.
Let the eigenvalues of $A$ and $B$ be denoted by $\{\lambda_i\}$ and $\{\hat{\lambda}_i\}$, respectively, and assume that they have been arranged in nonincreasing order $\lambda_1 \ge \cdots \ge \lambda_n$ and $\hat{\lambda}_1 \ge \cdots \ge \hat{\lambda}_{n-1}.$ Then
\begin{center}
$\lambda_i \ge \hat{\lambda}_i \ge \lambda_{i+1}$ for $ i=1,2,\ldots, n-1$.
\end{center}
\end{theorem}
A simple consequence of the Interlacing Theorem is the following lemma:
\begin{lemma} \label{lemm:upsize} Let $G$ be a graph on $n$ vertices. Let $\hat{G}$ be a graph formed by adding a vertex to $G$ and any number of edges between the new vertex and any vertices of $ G$. Let $A$ and $\hat{A}$ be the adjacency matrices of $G$ and $\hat{G}$ respectively. Then $i_+(\hat{A})+i_0(\hat{A}) \ge i_+(A)+i_0(A)$ and $i_-(\hat{A})+i_0(\hat{A}) \ge i_-(A)+i_0(A)$.
\end{lemma}

The following result on the sum of an Hermitian and a positive semidefinite matrix is from \cite{CJ}.
\begin{lemma} \label{lemm:possemidef} Let $A,B \in M_n$ be Hermitian. Assume that $B$ is positive semidefinite and that the eigenvalues of $A$ and $A+B$ are arranged in non-increasing order. Then
$$ \lambda_k(A)+\lambda_n(B)  \le \lambda_k(A+B) \text{ for } k=1,2, \ldots ,n.$$

\end{lemma}
Finally, we will need certain Ramsey numbers. The Ramsey number $R(m,n)$ is the minimum number of vertices such that all graphs of order $R(m,n)$ or more have either an independent set of size $m$ or a complete graph of order $n$ as an induced subgraph. Ramsey numbers are known to exist for all $(m,n)$, however the exact values are not known beyond $R(3,9)$ and $R(4,5)$ (\cite{Ex,GG,MR,Rams}).

\section{The Existence of Bounds}

Our primary question may be stated as follows: Let $k$ be a given positive integer. Is there an integer $n$ for which the adjacency matrix of any graph of order at least $n$ has at least $k$ nonpositive eigenvalues? We shall see that such an $n$ exists for each $k$. We denote this minimum $n$ by $NPO(k)$. For smaller numbers of vertices, some graphs have fewer than $k$ nonpositive eigenvalues. Using the concept of inertia, we may give an alternative description for $NPO(k)$: Let $k$ be a given positive integer. Is there an integer $n$ for which the adjacency matrix $A$ of any graph on at least $n$ vertices satisfies $i_-(A)+i_0(A) \geq k$? This minimum size is just $NPO(k)$.\\

We start by proving that $NPO(k)$ exists for any $k$.
\begin{theorem} \label{thm:rams} We have $NPO(k) \leq R(k,k+1)$.
\end{theorem}
\begin{proof}

By the definition of Ramsey numbers, any graph of order $R(k,k+1)$ or greater has either an independent set of size $k$ or a complete graph of order $k+1$ as an induced subgraph.
\\
If $G$ has an independent set of size $k$, the $k-\text{by}-k$ zero matrix is a principle submatrix of $A(G)$. Using the Interlacing Theorem, it follows that $A(G)$ has at least $k$ nonpositive eigenvalues.
\\
If $G$ has the complete graph on $k+1$ vertices as an induced subgraph, the matrix $J_{k+1}-I_{k+1}$ is a principle submatrix of $A(G)$ ($J$ is the matrix all of whose entries are 1). The eigenvalues of $J_{k+1}-I_{k+1}$ are $-1$ of multiplicity $k$ and $k$ of multiplicity $1$. Therefore, from the Interlacing Theorem, $A(G)$ has at least $k$ eigenvalues (counting multiplicities) that are smaller than or equal to $-1$. In particular, $A(G)$ has at least $k$ nonpositive eigenvalues.
\end{proof}
In the following three corollaries, we use Ramsey numbers whose values have been determined in \cite{GG}, \cite{Ex,MR} and \cite{Ex} respectively.
\begin{corollary}
$NPO(3) \leq 9=R(3,4)$
\end{corollary}
\begin{corollary}
$NPO(4) \leq 25=R(4,5)$
\end{corollary}
\begin{corollary}
$NPO(5) \leq R(5,6) \leq 87$
\end{corollary}
Since the Ramsey numbers are known to exist and be finite for all parameters (as shown in \cite{Rams}), Theorem ~\ref{thm:rams} demonstrates that for any positive integer $k$, $NPO(k)$ exists. This result is limited, however, by the rapid increase in the Ramsey numbers. While it demonstrates that such a bound exists, it is far from the actual value.\\

We next use generalized Ramsey numbers to get sharper bounds and start with the following lemma(\cite{Smith}):
\begin{lemma} \label{lemm:for k4first} Let $G$ be a graph on $n$ vertices that is the complement of a disjoint union of any number of complete graphs. Then $G$ has at least $n-1$ nonpositive eigenvalues.
\end{lemma}
This lemma has several useful consequences, but we will need the concept of generalized Ramsey numbers. Instead of using $R(k,k+1)$, we may obtain an upper bound on $NPO(k)$, by taking the smallest number $r$ for which all the graphs of order $r$ contain as an induced subgraph, at least one of the graphs in $S$, with $S$ a set of graphs that contains the complete graph of order $k+1$, the empty graph (i.e. an independent set) of order $k$, and other graphs that have $k$ nonpositive eigenvalues. Such $r$ would be also an upper bound on $NPO(k)$, and in many cases it may be much better than $R(k,k+1)$. In order to use this concept, we start with a more general definition of Ramsey numbers, which can be found in \cite{small Ramsey}: Let $G$ and $H$ be two graphs. The generalized Ramsey number $R(G,H)$ is the minimum number of vertices such that all graphs of order at least $R(G,H)$ have either a subgraph that is isomorphic to $G$ or the complement has a subgraph that is isomorphic to $H$. Note that in both cases, these subgraphs \textbf{are not necessarily induced subgraphs}.
Using this new definition, lemma ~\ref{lemm:for k4first}, and the values of generalized Ramsey numbers which can be found in \cite{small Ramsey}, we obtain the following bounds that are better than the previous ones. We denote by $K_n \setminus e$ the graph $K_n$ after removing one edge:

\begin{corollary}
$NPO(3) \leq 7=R(K_4 \setminus e,K_3)$
\end{corollary}
\begin{corollary}
$NPO(4) \leq 19=R(K_5 \setminus e,K_4)$
\end{corollary}
\begin{corollary}
$NPO(5) \leq R(K_6 \setminus e,K_5) \leq 67$
\end{corollary}

We also mention another lemma and theorem from \cite{small Ramsey}, that will help us later calculate $NPO(3)$ and $NPO(5)$:
\begin{lemma} \label{lemm:special ramsey}
$R(K_{2,2},K_{1,3})=6$
\end{lemma}
\begin{theorem} \label{thm:special ramsey2}
$R(K_4 \setminus e,K_5)=16$
\end{theorem}

We now improve the prior values for certain numbers of nonpositive eigenvalues.
\begin{lemma} $NPO(1)=1$
\end{lemma}
\begin{proof}
Let $G$ be a graph which has at least one vertex. $A(G)$ is hollow (i.e. has 0 diagonal) by the definition of an adjacency matrix. Thus, the 1-by-1 principal submatrix of $A(G)$ is $\textbf{0}$. Therefore, from the Interlacing theorem, $A(G)$ has at least $1$ nonpositive eigenvalue.
\end{proof}
\begin{lemma}\label{lemm:npotwothree} $NPO(2)=3$
\end{lemma}
\begin{proof}
We apply Theorem~\ref{thm:rams}, using the fact that $R(2,3)=3$ to get that $NPO(2) \leq 3$. In addition, $K_2$ has only one nonpositive eigenvalue, and therefore $NPO(2)=3$.
\end{proof}
While the bounds for one and two nonpositive eigenvalues are simple to determine, complexity increases dramatically beyond this point. But it does motivate the attempt to determine a bound better than that given by Ramsey numbers.

\section{Exact values for $k=3$, $k=4$, $k=5$}

We start with the following two lemmas:
\begin{lemma} \label{lemm:build1} Let $G$ be a graph on $m$ vertices such that the adjacency matrix $A(G)$ has at least $k$ nonpositive eigenvalues. Then, any graph $H$ on $p > m$ vertices, that has $G$ as an induced subgraph, has at least $k$ nonpositive eigenvalues.
\end{lemma}
\begin{proof}
This is a simple consequence of Lemma~\ref{lemm:upsize}.
\end{proof}
\begin{corollary} \label{coro:build2} Let $n$ be an integer such that for any graph $G$ on $n$ vertices, the adjacency matrix $A(G)$ has at least $k$ nonpositive eigenvalues. Then for any graph $H$ on $m \ge n$ vertices, the adjacency matrix $A(H)$ has at least $k$ nonpositive eigenvalues.
\end{corollary}
\begin{proof}
There exists a graph $G$ on $n$ vertices such that $H$ has $G$ as an induced subgraph. Since $G$ has at least $k$ nonpositive eigenvalues, by Lemma~\ref{lemm:build1}, $H$ has at least $k$ nonpositive eigenvalues.
\end{proof}
The key consequence of these two lemmas is the fact that once a number $n$ is found such that any graph on $n$ vertices has an adjacency matrix with at least $k$ nonpositive eigenvalues, all graphs on more vertices do, as well. That integer $n$ will bound $NPO(k)$ from above. If we find an example of a graph on $n-1$ vertices whose adjacency matrix has less than $k$ nonpositive eigenvalues, we may conclude that $NPO(k)=n$.\\

We may now determine the exact value of $NPO(3)$:
\begin{theorem} \label{thm:k3} $NPO(3)=6$
\end{theorem}

\begin{proof}
First, since $C_5$ has only 2 nonpositive eigenvalues, $NPO(3)>5$. Using corollary ~\ref{coro:build2}, It is enough to show that for any graph $G$ on 6 vertices, $A(G)$ has at least 3 nonpositive eigenvalues. Let $G$ be a graph of order 6. Using Lemma ~\ref{lemm:special ramsey}, either $G$ has a subgraph that is isomorphic to $K_{2,2}$, or $G$ has an induced subgraph that is isomorphic to $K_1 \cup H$, in which $H$ is a graph of order 3. In the first case, using Lemma ~\ref{lemm:for k4first} we get that $G$ has an induced subgraph of order 4 with 3 nonpositive eigenvalues, and we are done. In the second case, since $NPO(2)=3$, $H$ has at least 2 nonpositive eigenvalues, and hence $K_1 \cup H$ has at least three nonpositive eigenvalues, and again, we are done.
\end{proof}

Note that the actual value of $NPO(3)$ is much lower than the Ramsey bound. Before we continue with the value of $NPO(4)$, we need the following lemma:
\begin{lemma} \label{lemm:for k4second} A graph $G$ on 4 vertices with just 2 nonpositive eigenvalues must satisfy $\delta(G)=1$. The only graph $G$ on 5 vertices with just 2 nonpositive eigenvalues is $C_5$.
\end{lemma}
\begin{proof}
First of all, notice that a graph on four or five vertices must have at least two nonpositive eigenvalues, since $NPO(2)=3$. Now, Let $G$ be a graph on 4 vertices. If $\delta(G) \geq 2$ then by Lemma~\ref{lemm:for k4first}, $G$ has 3 nonpositive eigenvalues. If $\delta(G) = 0$, then since $NPO(2)=3$, and because of the Interlacing Theorem we get that $G$ has 3 nonpositive eigenvalues. Therefore, a graph $G$ on 4 vertices with just 2 nonpositive eigenvalues must satisfy $\delta(G)=1$. Now consider the second statement. Let $G$ be a graph on 5 vertices. Since $NPO(2)=3$, it follows that $\delta(G) \geq 2$ (otherwise, by looking at the subgraph induced by the vertex with degree 1, and the three vertices that are not connected to it we get that $G$ must have at least 3 nonpositive eigenvalues). If the degrees of all the vertices are 2, then $G=C_5$ and we are done. Otherwise, there exists a vertex of degree at least 3, and using Lemma ~\ref{lemm:for k4first} and that $K_{1,3}, K_4 \setminus e$ and $K_4$ cannot be induced subgraphs (since otherwise we are done) we get that $G$ contains an induced subgraph on 4 vertices that is $K_3$ with a pedant vertex. Since $\delta(G) \geq 2$, this pedant vertex has to be connected to the remaining vertex of $G$, and from here it is easy to check that for all the possible connections to the remaining vertex, we get that $G$ has at least 3 nonpositive eigenvalues. Hence, a graph on 5 vertices with just 2 nonpositive eigenvalues must be $C_5$.
\end{proof}
We may now determine the value of $NPO(4)$:
\begin{theorem} \label{thm:k4} $NPO(4)=10$
\end{theorem}
\begin{proof}
The proof is composed of two parts: $NPO(4) \leq 10$, and $NPO(4) > 9$. We start with the first part. Using Corollary ~\ref{coro:build2}, it is enough to show that for any graph $G$ on 10 vertices, $A(G)$ has at least 4 nonpositive eigenvalues. Suppose in contradiction that there exists a graph $G$ on 10 vertices with less than 4 nonpositive eigenvalues. We divide the proof into several cases.
\begin{enumerate}
  \item $\delta(G) \leq 3$
  \item $\delta(G) \geq 6$
  \item $\delta(G) = 5$
  \item $\delta(G) = 4$
\end{enumerate}
In case 1, Let us look at a vertex $v$ in $G$ whose degree is at most 3. There are at least 6 vertices that are not connected to $v$. Let $\hat{G}$ be a subgraph of $G$ induced by these vertices. From Theorem ~\ref{thm:k3}, $A(\hat{G})$ has at least 3 nonpositive eigenvalues. By the definition of $v$, $\hat{G} \cup v$ is also an induced subgraph of $G$, and since $v$ is not connected to any vertex of $\hat{G}$, the adjacency matrix of $\hat{G} \cup v$ has at least 4 nonpositive eigenvalues, and therefore by the Interlacing Theorem $G$ has at least 4 nonpositive eigenvalues, which contradicts the assumption.\\
In case 2, let us look at a vertex $v$ in $G$ whose degree is $\delta(G)$, and recall that $\delta(G) \geq 6$. If $\delta(G) = 9$ then we are done, since $G$ has to be $K_{10}$, and hence its adjacency matrix has 9 nonpositive eigenvalues. So assume that $\delta(G) < 9$. In this case, there exists a vertex $u$ in $G$ that is not connected to $v$. There are at most 3 vertices (including $u$) that are not connected to $v$, and since the degree of $u$ is at least 6, $u$ and $v$ have at least 4 common neighbors. Therefore, $G$ has a subgraph (not necessarily induced) isomorphic to the graph in Figure ~\ref{first-fig}.
\begin{figure}[h!]
\caption{}
\centering
\includegraphics[width=0.3\textwidth]{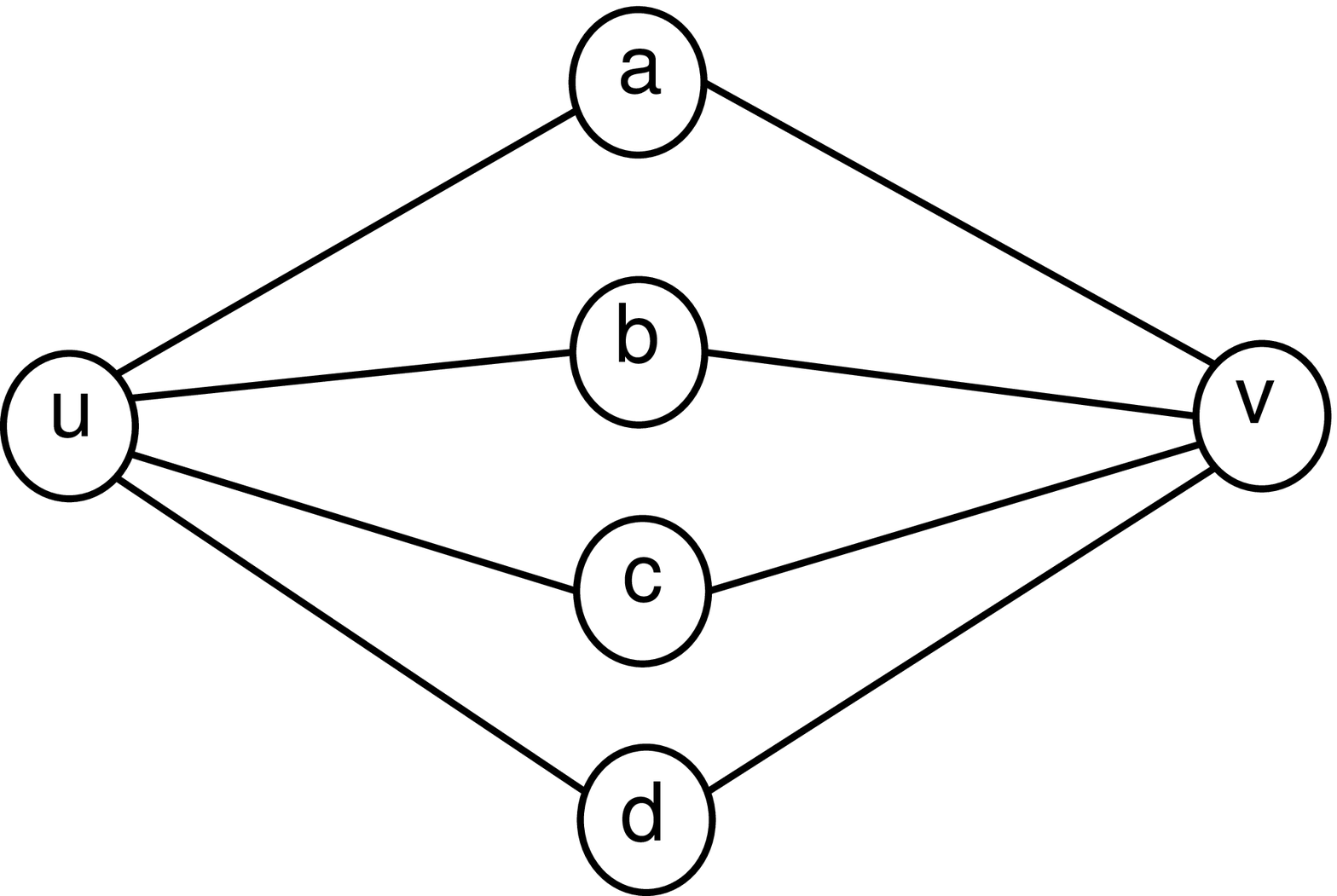}
\vspace{-0.85\baselineskip}
\label{first-fig}
\end{figure}
\\
Since the complete bipartite graph on $c$ vertices has $c-1$ nonpositive eigenvalues (which is a special case of Lemma ~\ref{lemm:for k4first}), and since $G$ has less than 4 nonpositive eigenvalues, $G$ has at least two of the following edges:\\ $\big\{a,b\big\},\big\{a,c\big\},\big\{a,d\big\},\big\{b,c\big\},\big\{b,d\big\},\big\{c,d\big\}$. If it has at least three of them, then either $\overline{K_2 \cup 3K_1}$  or $\overline{2K_2 \cup K_1}$ is an induced subgraph of $G$, and then we get a contradiction by Lemma ~\ref{lemm:for k4first}. If $G$ has exactly two edges among the 6 that were described above, then without loss of generality, either $\overline{2K_2 \cup K_1}$ or the graph in Figure~\ref{second-fig} is an induced subgraph of $G$.\\
\begin{figure}[h!]
\caption{}
\centering
\includegraphics[width=0.3\textwidth]{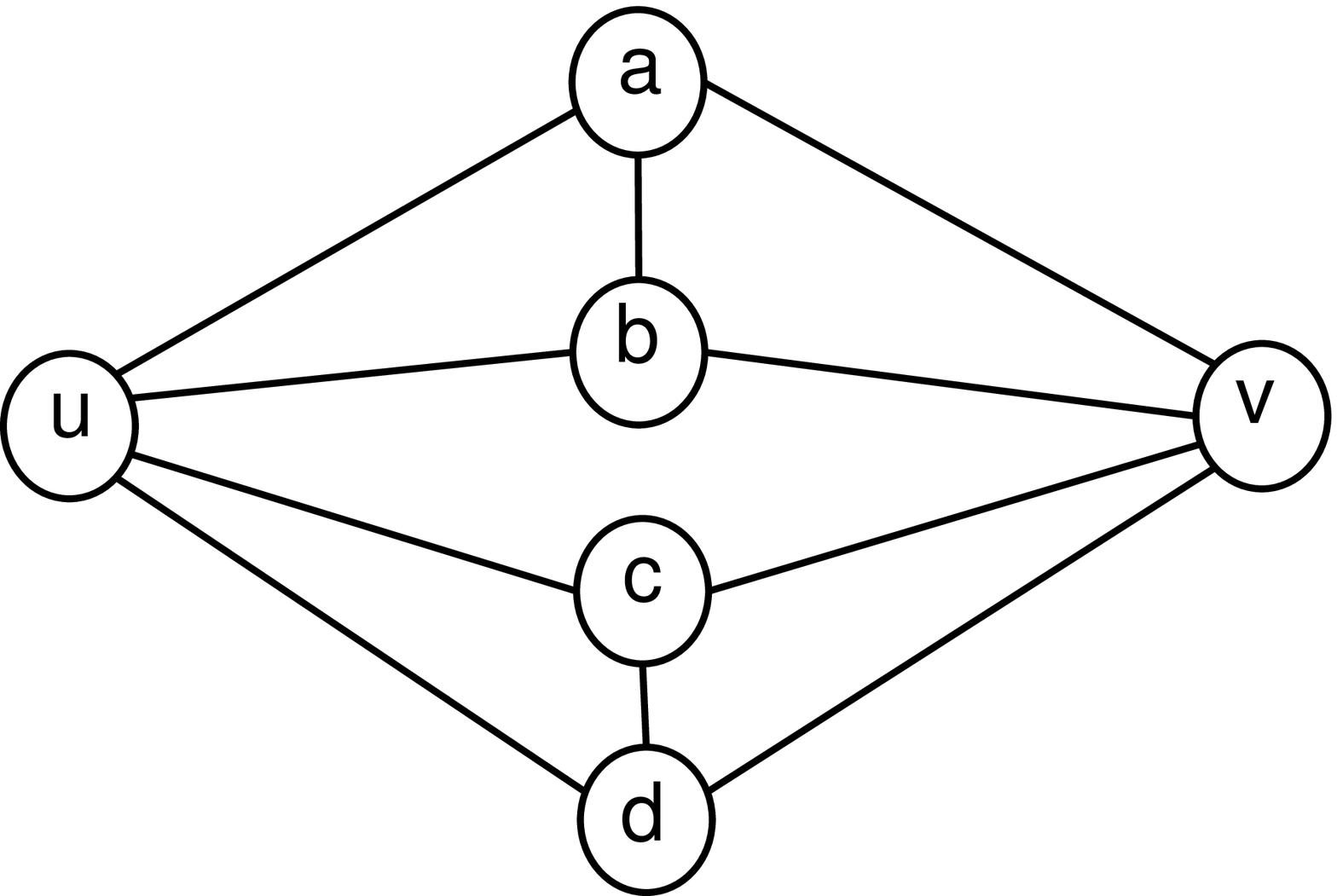}
\vspace{-0.85\baselineskip}
\label{second-fig}
\end{figure}
In the first case we get a contradiction. The eigenvalues of the graph in Figure~\ref{second-fig} are $\big\{\frac{1}{2} (1 - \sqrt{33}), -1, -1, 0,1, \frac{1}{2} (1 + \sqrt{33})\big\}$, therefore from the Interlacing Theorem $G$ has at least 4 nonpositive eigenvalues, which again contradicts the assumption, and hence case 2 is impossible.\\
In case 3, let us look at a vertex $v$ in $G$ whose degree is $5$. There are 4 vertices that are not connected to $v$. Let $\hat{G}$ be a subgraph of $G$ induced by these vertices. If $\hat{G}$ has at least 3 nonpositive eigenvalues, then we can continue in the same way as in case 1 and get a contradiction. Therefore $\hat{G}$ has to have at most 2 nonpositive eigenvalues. From Lemmas ~\ref{lemm:npotwothree} and ~\ref{lemm:for k4second}, $\hat{G}$ has a vertex with degree one. we denote this vertex by $u$. Since the degree of $u$ in $G$ is at least 5, $u$ and $v$ have at least 4 common neighbors. From here we continue in the same way as in case 2, and we get a contradiction.\\
The last case is case 4. Let $v$ be a vertex of degree 4. There are 5 vertices that are not connected to $v$. Let $\hat{G}$ a subgraph of $G$ induced by these vertices. For the same reason as before, $\hat{G}$ has to have at most 2 nonpositive eigenvalues (otherwise we get a contradiction). Therefore, from Lemma~\ref{lemm:for k4second}, $\hat{G}$ is $C_5$. Let $u$ be some vertex in $\hat{G}$. Since $\delta(G) = 4$, $u$ and $v$ have at least 2 common neighbors. Therefore, $G$ has a subgraph (not necessarily induced) that is isomorphic to the one in Figure ~\ref{third-fig}.\\
\begin{figure}[h!]
\caption{}
\centering
\includegraphics[width=0.6\textwidth]{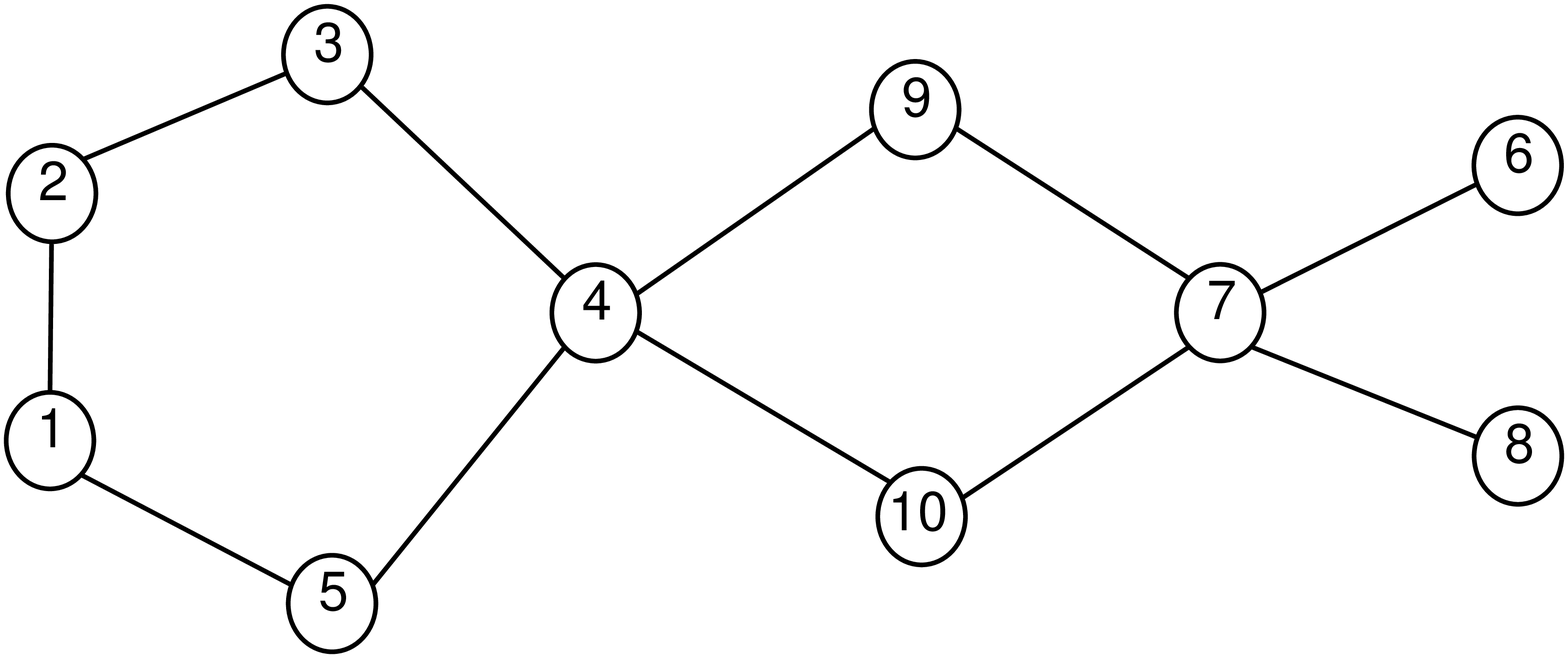}
\vspace{-0.85\baselineskip}
\label{third-fig}
\end{figure}
Let us look at vertices 1,2,4,7,9,10. Since vertices $1,2,4$ are part of $\hat{G}$, which is an induced subgraph, and since the degree of vertex 7 is 4, the only possible edges among the set of vertices $1,2,4,7,9,10$ are\\ $\big\{1,9\big\},\big\{1,10\big\},\big\{2,9\big\},\big\{2,10\big\},\big\{9,10\big\}$. If at least one of vertices $1,2$ is connected to both $9$ and $10$, then the subgraph induced by this vertex and vertices $4,7,9,10$ has 4 nonpositive eigenvalues by Lemma~\ref{lemm:for k4first}, which leads us to a contradiction. On the other hand, since in all the possible cases the subgraph induced by vertices $4,7,9,10$ has 3 nonpositive eigenvalues, if at least one of vertices $1,2$ is not connected neither to $9$ nor to $10$, then the subgraph induced by $1,2,4,7,9,10$ has at least 4 nonpositive eigenvalues and we get a contradiction. Therefore each of vertices 1 and 2 is connected to exactly one of the vertices 9 and 10. Hence, there are four options up to isomorphism for the subgraph induced by $1,2,4,7,9,10$. It is easy to check that only in one of them there are less than 4 nonpositive eigenvalues. Therefore, up to isomorphism, $G$ has a subgraph that is isomorphic to the one in Figure ~\ref{four-fig}.\\
\begin{figure}[h!]
\caption{}
\centering
\includegraphics[width=0.6\textwidth]{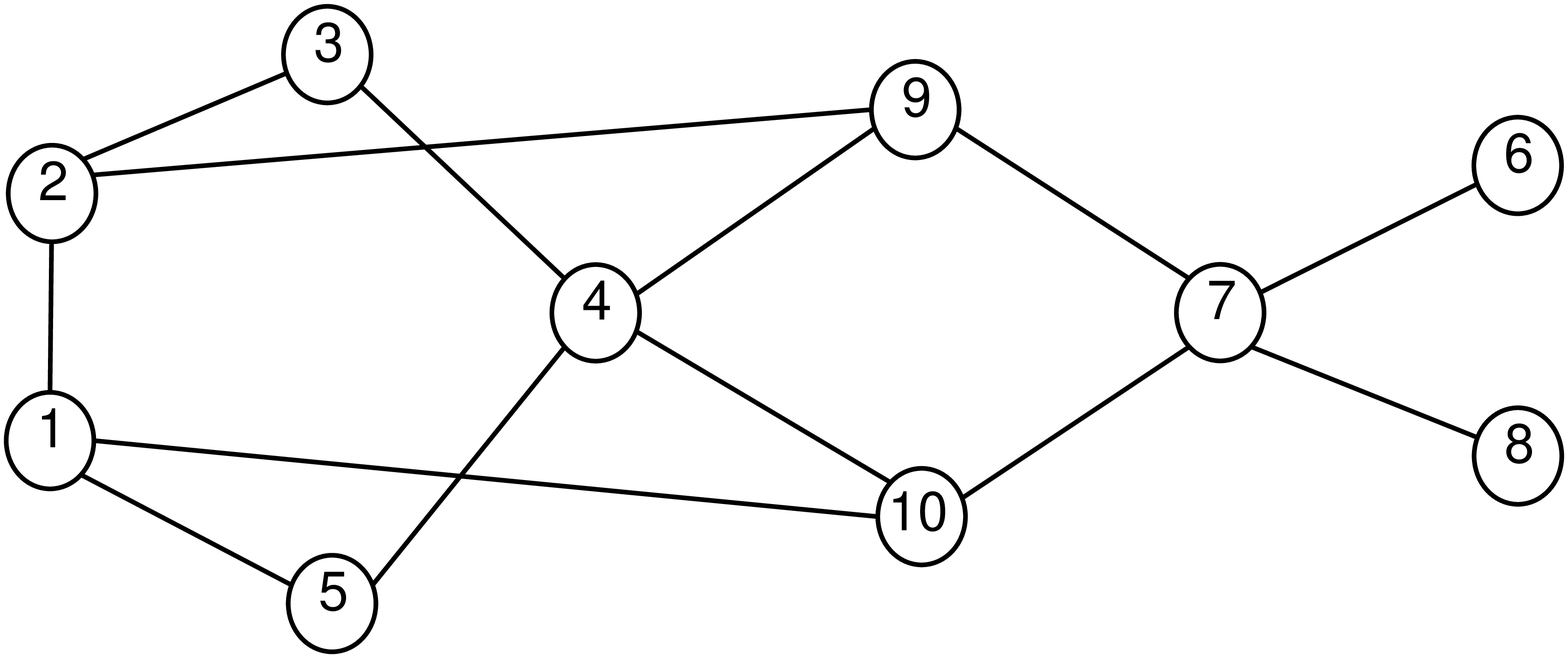}
\vspace{-0.85\baselineskip}
\label{four-fig}
\end{figure}
Let us look now at the subgraph induced by all the vertices except 6 and 8. There are only four optional edges that may be added (other than those already accounted for): $\big\{3,9\big\},\big\{3,10\big\},\big\{5,9\big\},\big\{5,10\big\}$. If all of them are added, or none of them is added, then we are done by looking at the induced subgraph on vertices $3,4,5,9,10$, and using Lemma~\ref{lemm:for k4first}. Therefore, up to isomorphism, there are only 8 cases that we didn't check yet for the subgraph of $G$ that is induced by all the vertices except 6 and 8. It is easy to check that only 3 of them have less than 4 nonpositive eigenvalues. These cases are illustrated in Figures ~\ref{five-fig},~\ref{six-fig} and ~\ref{seven-fig}.

\begin{figure}[h!]
\caption{}
\centering
\includegraphics[width=0.6\textwidth]{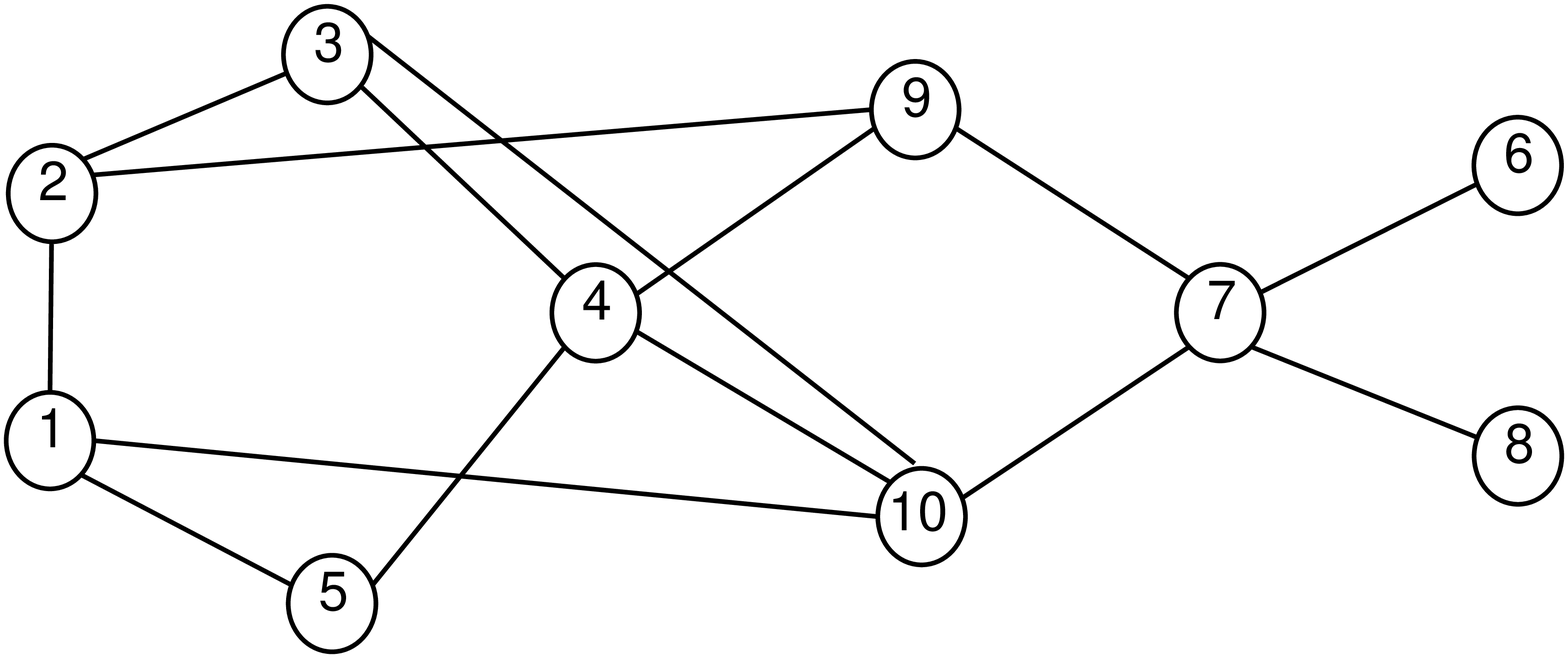}
\vspace{-0.2\baselineskip}
\label{five-fig}
\end{figure}
\begin{figure}[h!]
\caption{}
\centering
\includegraphics[width=0.6\textwidth]{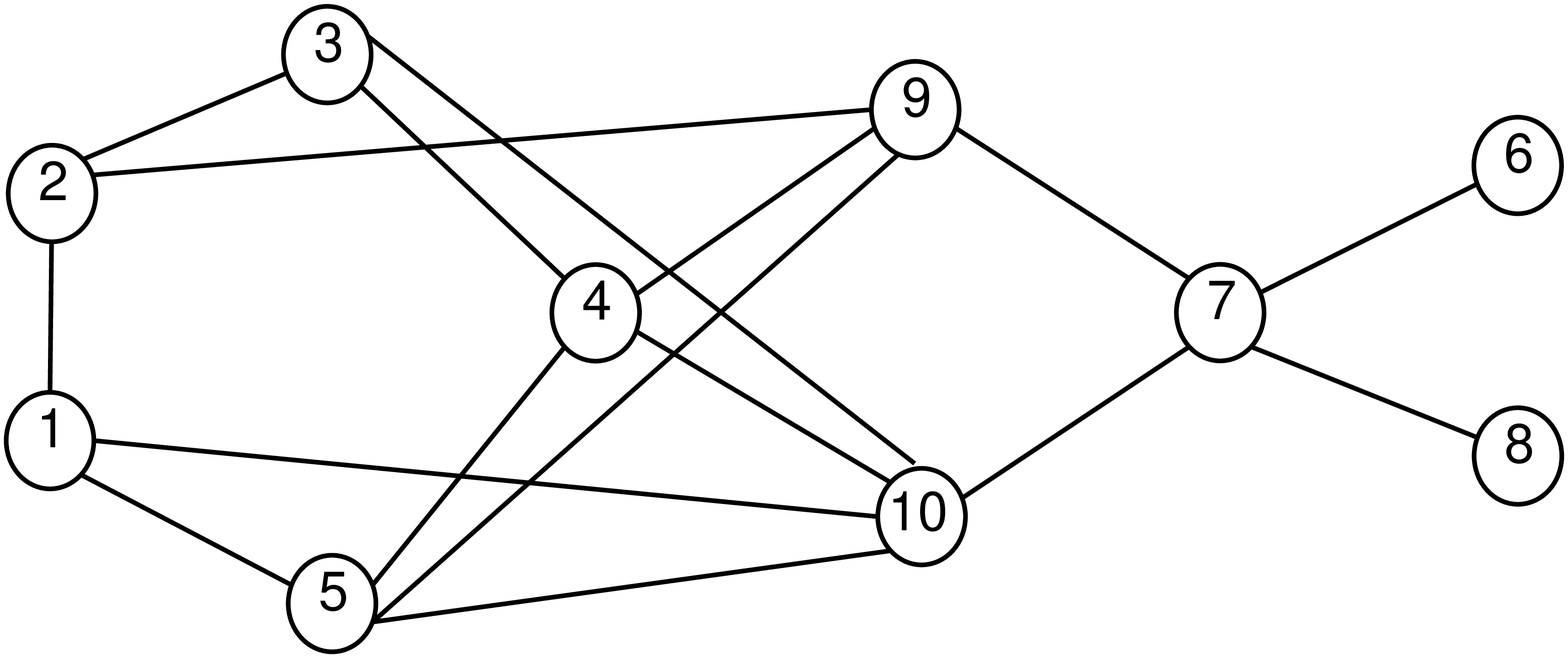}
\vspace{-1\baselineskip}
\label{six-fig}
\end{figure}
\begin{figure}[h!]
\caption{}
\centering
\includegraphics[width=0.6\textwidth]{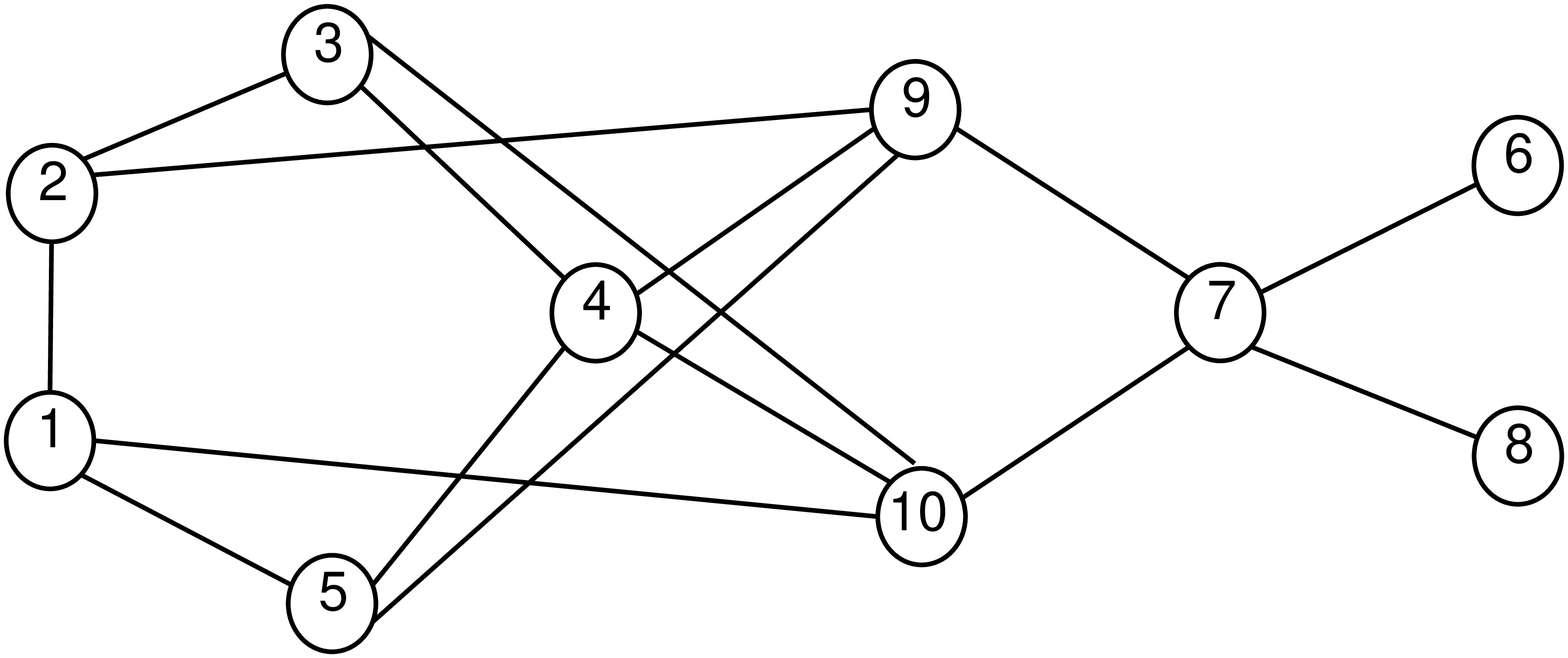}
\vspace{-1\baselineskip}
\label{seven-fig}
\end{figure}

We start by examining Figure~\ref{five-fig}. Since $\delta(G) = 4$, vertex 5 is connected both to vertices 6 and 8. In addition, vertex 3 has to be connected to at least one of vertices 6 and 8. If it is connected to only one of them, then from Lemma~\ref{lemm:for k4second}, either the induced subgraph on vertices $1,5,9,7,8$, or the induced subgraph on vertices $1,5,9,7,6$, has to be a cycle. Both of these cases are impossible since vertices 1 and 9 are not connected. Therefore, vertex 3 is connected both to vertices 6 and 8. Let us look at the subgraph induced by vertices $3,5,7,6,8$. Using Lemma~\ref{lemm:for k4first}, this subgraph has 4 nonpositive eigenvalues (in both possible cases where either there is or there is no edge between vertices 6 and 8), and we get a contradiction. Therefore the situation illustrated in Figure~\ref{five-fig} is impossible.\\
The next case is Figure~\ref{six-fig}. Since $\delta(G) = 4$, vertex 2 is connected to at least one of vertices 6 and 8. Without loss of generality, we may assume that vertex 2 is connected to vertex 6. If the degree of vertex 2 is four, then using Lemma ~\ref{lemm:for k4second}, we get that the induced subgraph on vertices $4,5,10,7,8$ must be a cycle, which is not the case, so we get a contradiction. Therefore, vertex 2 is also connected to vertex 8. Now, let us look at vertex 4. If this vertex is not connected to 6 nor to 8, then in all the possible cases, the subgraph induced by vertices $2,4,7,6,8$ has 4 nonpositive eigenvalues, and we get a contradiction. Therefore, without loss of generality, we may assume that vertex 4 is connected to vertex 6. Finally, Using Lemma ~\ref{lemm:for k4first}, we get that the subgraph induced by vertices $6,9,4,2,7$ has 4 nonpositive eigenvalues in all the possible cases (there are two cases, either vertices 6 and 9 are connected, or not), therefore the situation that is illustrated in Figure~\ref{six-fig} is impossible.\\
The last case that we have is the situation illustrated in Figure~\ref{seven-fig}. Since $\delta(G) = 4$, vertex 5 is connected to at least one of vertices 6 and 8. Without loss of generality, we may assume that vertex 5 is connected to vertex 8. In addition, vertex 3 is connected to at least one of vertices 6 and 8. If it is connected only to vertex 6, then from Lemma~\ref{lemm:for k4second}, the subgraph induced by vertices $1,5,9,7,8$ has to be a cycle, which is not the case. Therefore, vertex 3 has to be connected to vertex 8. Now we have two cases: Either vertex 3 is connected to vertex 6, or not. If they are connected, then vertex 5 has to be connected to vertex 6 (otherwise, the subgraph induced by vertices $2,3,10,7,6$ has to be a cycle, which is not the case). Now, let us look at the subgraph induced by vertices $3,5,7,6,8$. Using Lemma ~\ref{lemm:for k4first}, we get that this subgraph has 4 nonpositive eigenvalues in all the possible cases, which leads us to a contradiction. Therefore, vertices 3 and 6 cannot be connected, and the degree of vertex 3 is four. From Lemma~\ref{lemm:for k4second}, the subgraph induced by vertices $1,5,9,7,6$ has to be a cycle, and hence vertices 1 and 6 are connected, and in addition there is no edge between vertices 5 and 6 and between vertices 6 and 9. Hence, the degree of vertex 5 is four, so the subgraph induced by vertices $2,3,7,10,6$ has to be a cycle. Therefore vertices 2 and 6 are connected, and there is no edge between vertices 6 and 10. Let us look at the subgraph induced by vertices $6,10,1,7,4$. If vertices 4 and 6 are connected, then using Lemma~\ref{lemm:for k4first} we get that this subgraph has 4 nonpositive eigenvalues, which leads us to a contradiction. Hence vertices 4 and 6 are not connected. Therefore the degree of vertex 6 is four, and it must be connected to vertex 8. So by Lemma~\ref{lemm:for k4second}, the subgraph induced by vertices $3,4,5,9,10$ has to be a cycle, which is not the case, so we get a contradiction. Therefore case 4 is impossible. So, in conclusion, after checking all the possible cases, we get that any graph on 10 vertices has to have at least 4 nonpositive eigenvalues. We conclude the proof by giving two examples of graphs of order 9 whose adjacency matrices have only 3 nonpositive eigenvalues, which means that $NPO(4)>9$. The examples are given in Figure ~\ref{fig:nineac}.
The eigenvalues of $A(G_2)$ are
\begin{center}
$\big\{-2.4142,-2.4142,-2.0000,0.4142,0.4142,0.5858,1.0000,1.0000,3.4142\big\},$
\end{center}
and the eigenvalues of $A(G_3)$ are
\begin{center}
$\big\{-2.4142,-2.4142,-2.1413,0.4142,0.4142,0.5151,1.0000,1.0000,3.6262\big\}.$
\end{center}
\begin{figure}[H]
\caption{}
\begin{center}
\mbox{
    \leavevmode
	\subfigure [{$G_2$}]
	{ 
	  \includegraphics[width=5cm]{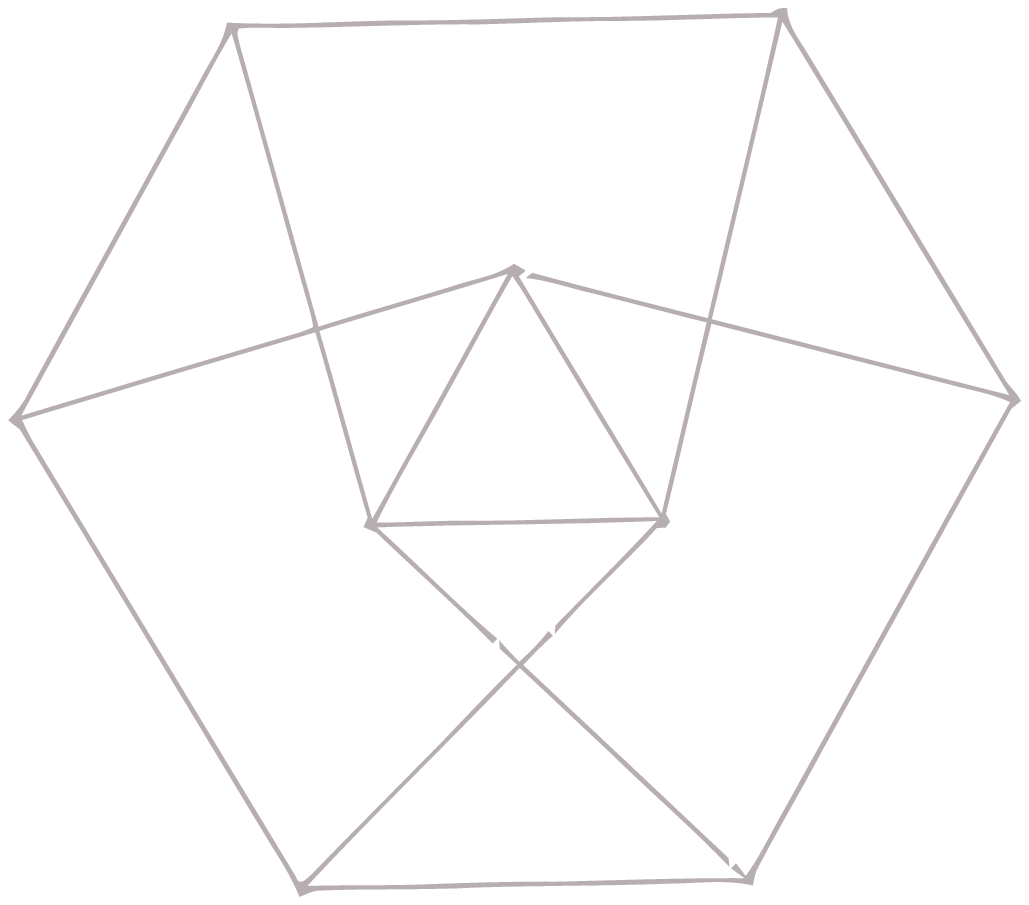} }

	\leavevmode
	\subfigure [{$G_3$}]
	{ 
	  \includegraphics[width=5cm]{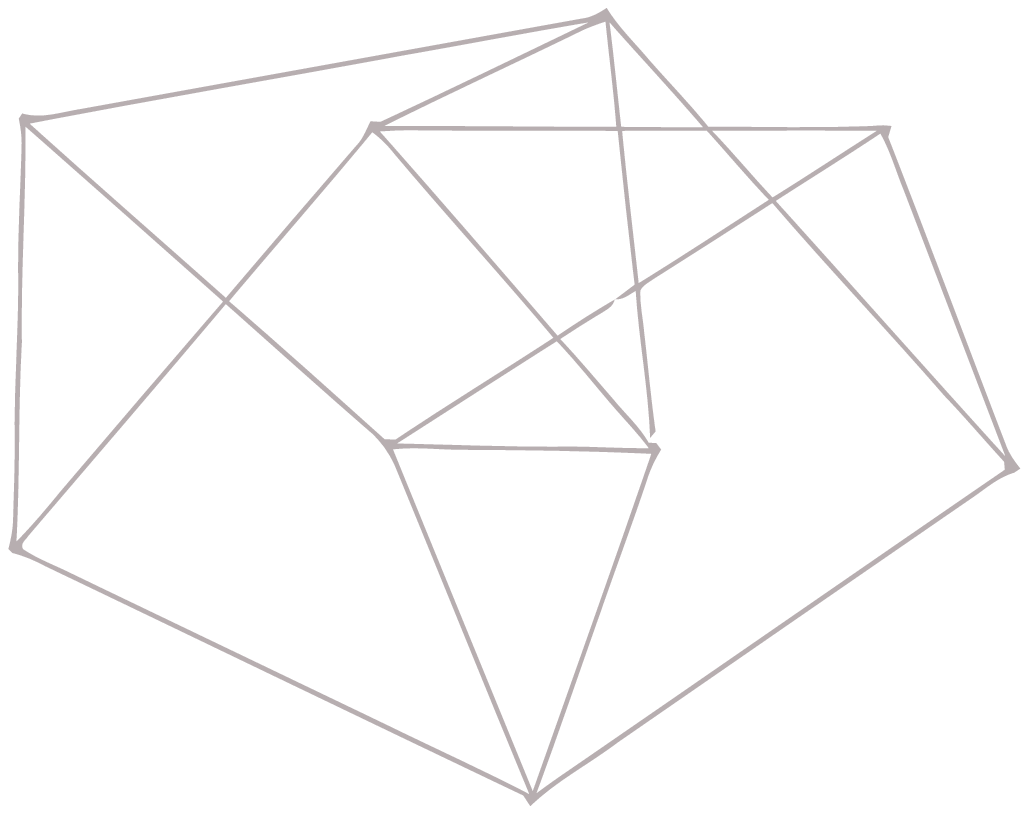} }
\label{fig:nineac}}
\end{center}
\end{figure}
\end{proof}
Our next goal is to determine the exact value of $NPO(5)$. We start with the following Lemma:
\begin{lemma} \label{lemm:k5} $NPO(5)>15$
\end{lemma}
\begin{proof}
We prove this by giving two graphs, each of which has 15 vertices and only 4 nonpositive eigenvalues. Both have very interesting structure. In addition, structure similar to that of the first graph will be presented in section 5 as a part of a more general bound for $NPO(k)$. The adjacency matrices of graphs $G_4$ and $G_5$ (see Figure ~\ref{fig:fifteenbc}) , both on $15$ vertices , have $4$ nonpositive eigenvalues. The eigenvalues of $G_4$ are\\
\big\{-3.3028 ,-3.3028 ,-3.3028, -3.3028 ,0.3028 ,0.3028 ,0.3028 ,0.3028 ,0.6277 ,1 ,1 ,1 ,1 ,1 , 6.3723\big\},\\
and the eigenvalues of $G_5$ are\\
\big\{-3.3028 ,-3.3028 ,-3.3028 ,-3.3028 ,0.3028 ,0.3028 , 0.3028 ,0.3028 ,0.3542 ,1 ,1 ,1 ,1 ,2 , 5.6458\big\}.\\
\begin{figure}[H]
\caption{}
\begin{center}
\mbox{
    \leavevmode
	\subfigure [{$G_4$}]
	{ 
	  \includegraphics[width=6.5cm]{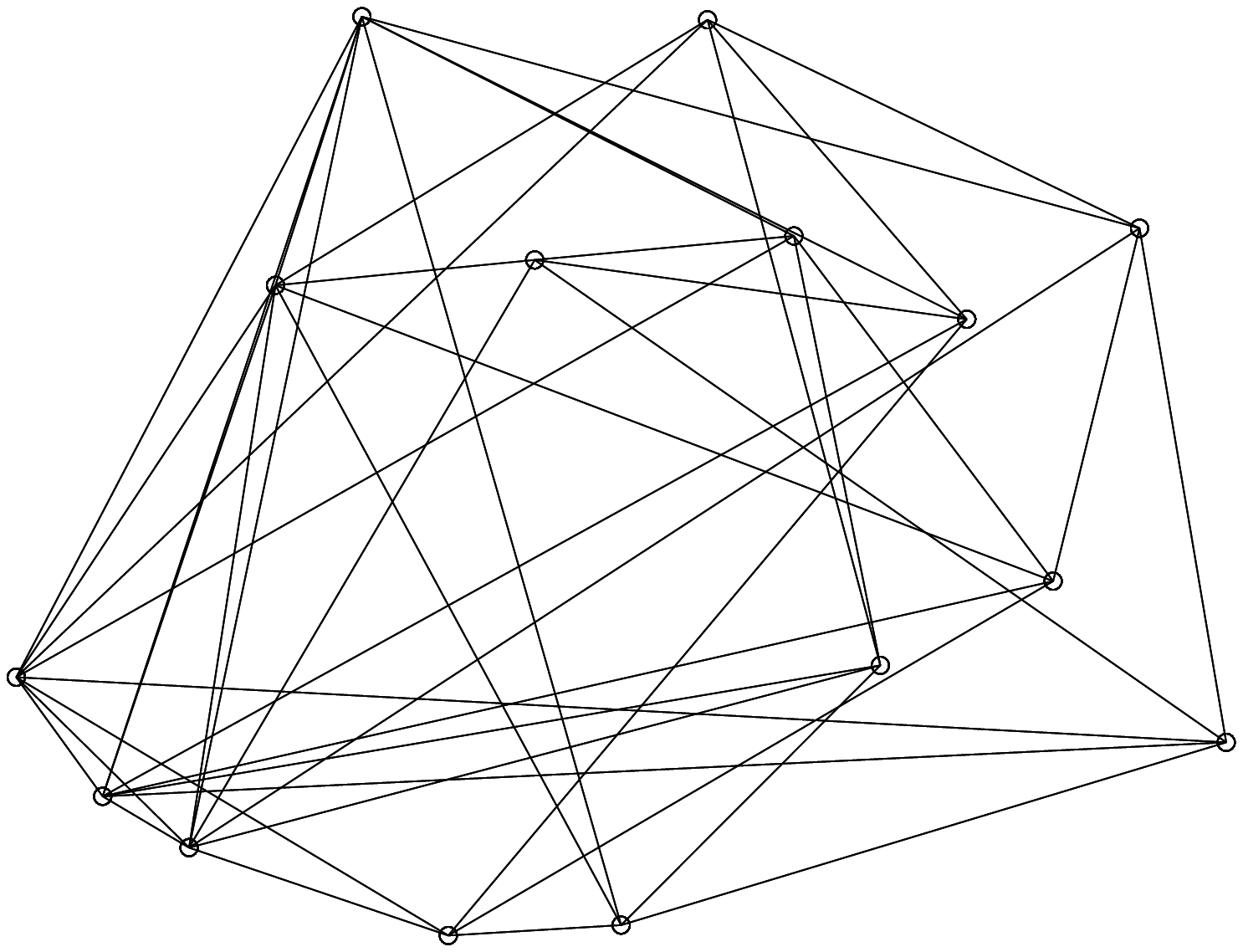} }

	\leavevmode
	\subfigure [{$G_5$}]
	{ 
	  \includegraphics[width=6.5cm]{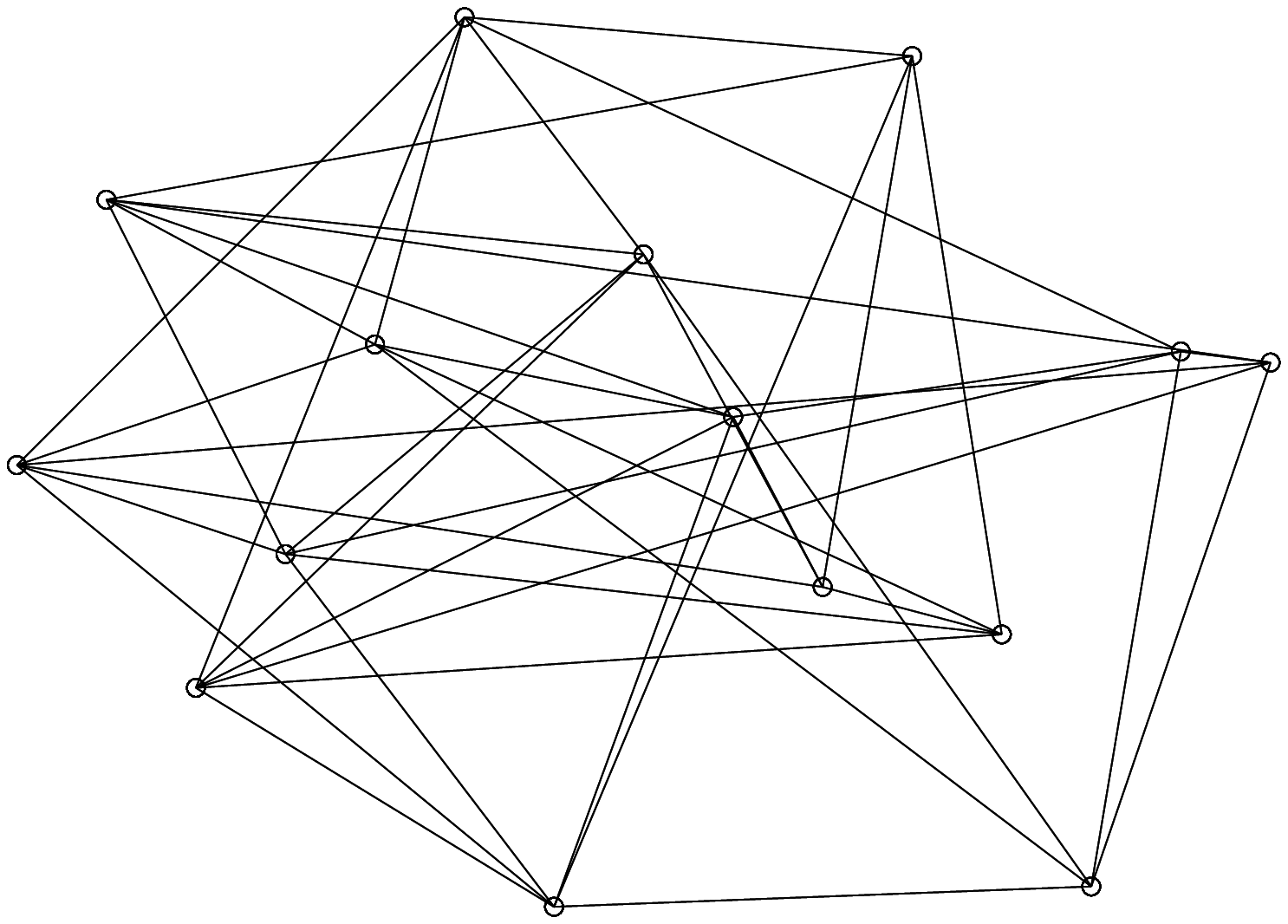} }
\label{fig:fifteenbc}
}
\end{center}
\end{figure}
 We describe now the structure of each, starting with $G_4$. First, Let $G$ be the disjoint union of the Petersen Graph and $K_5$, which has 15 vertices. The independence number of the Petersen Graph is four, and it has exactly five different maximal independent sets. We obtain $G_4$ by connecting each vertex in $K_5$ to all vertices of one of the five independent sets, such that each vertex in $K_5$ is connected to a different independent set.\\
In order to get $G_5$, we start with the Paley graph on 9 vertices, denoted by $P(9)$. We construct a graph $G$ by taking the disjoint union of $P(9)$ and $2K_3$. The independence number of $P(9)$ is 3, and it has six different maximal independent set. Moreover, it is possible to label the vertices of $P(9)$ with $1,2, \ldots, 9$ so that the six maximal independent sets are $\big\{1,5,9\big\},\big\{2,6,7\big\},\big\{3,4,8\big\},\big\{1,6,8\big\},\big\{2,4,9\big\},\big\{3,5,7\big\}$. Note that the first three sets are disjoint, and the last three sets are also disjoint. We obtain $G_5$ in the following way: take the first $K_3$ and connect each one of its vertices to an independent set, such that the first vertex would be connected to vertices $1,5,9$, the next vertex would be connected to vertices $2,6,7$, and the last to vertices $3,4,8$. Each vertex in the second $K_3$ would be connected to one of the three maximal independent sets that remain (each one of the three connected to a different independent set).

\end{proof}
Now we may determine the exact value of $NPO(5)$:
\begin{theorem} \label{thm:k5} $NPO(5)=16$
\end{theorem}
\begin{proof}
We already know that $NPO(5)>15$, so that it is enough to show that $NPO(5) \leq 16$. Let $G$ be a graph of order 16. First, since $NPO(4)=10$, if $\delta(G)<6$ then $G$ has to have at least 5 nonpositive eigenvalues. Therefore, we may assume that $\delta(G) \geq 6$. In addition, if the independence number of $G$ is 5 or more, then by the Interlacing Theorem, $G$ has to have at least 5 nonpositive eigenvalues. Therefore, we may also assume that the independence number of $G$ is smaller than 5. Using~\ref{thm:special ramsey2}, $G$ has a subgraph that is isomorphic to $K_4 \setminus e$. Together with the assumption that $\delta(G) \geq 6$ we get that $G$ has a subgraph of order 7 (not necessarily induced) that can be obtained from $K_{1,6}$ by choosing some vertex of degree one and connecting it to two other vertices. Let us denote this subgraph by $H$. From this point, the proof is completed computationally. We start by identifying the graphs of order $7$ that have a subgraph that is isomorphic to $H$, and whose adjacency matrix has fewer than $5$ nonpositive eigenvalues. Note that if there exists a graph on 16 vertices with fewer than 5 nonpositive eigenvalues, it has to have one of these graphs as an induced subgraph (otherwise by the Interlacing Theorem, the graph has at least 5 nonpositive eigenvalues). Our program identified $68$ non-isomorphic graphs on 7 vertices with fewer than 5 nonpositive eigenvalues and with $H$ as a subgraph. It continued by checking all the possible options for adding a vertex to each one of these 68 graphs (by all the possible options we mean also all the possible options for adding an edge between this vertex and the vertices of the induced subgraphs of order 7). Using this process, we got the set of all the graphs of order 8 with fewer than 5 nonpositive eigenvalues and with $H$ as a subgraph. We continued the process with the same idea, each time getting the set of all the options for graphs on $n$ vertices ($9 \leq n \leq 16$) with less than 5 nonpositive eigenvalues and with $H$ as a subgraph. For this purpose, we used the high-performance computing resources of the College of William and Mary. For graphs of order 16, the program found that there is no graph of order 16 with less than 5 nonpositive eigenvalues and with $H$ as a subgraph. Since all the graphs of order 16 and with fewer than 5 nonpositive eigenvalues have $H$ as a subgraph, we can conclude that all the graphs on $16$ vertices have adjacency matrices with at least $5$ nonpositive eigenvalues.
\end{proof}

These results are very promising. They demonstrate that the exact value that can be determined is much smaller than the bound from Theorem~\ref{thm:rams} would indicate.

\section{Lower bound for $NPO(k)$}
Another approach to finding bounds, is to find "extreme" graphs, i.e. families of graphs with increasing numbers of vertices and with a small number of nonpositive eigenvalues. We have identified a number of extreme graphs, including a construction that can generalize to any size, which we present here. These give what appear to be very strong lower bound for $NPO(k)$.
We start with the definition of triangular numbers. The triangular number $T_n$ is a number that can be represented as a triangular grid of points where the first row contains one point, and each subsequent row contains one more point than the previous one such that there are $n$ rows in total. Another way to define $T_n$ is to sum all the positive integers which are smaller than or equal to $n$.\\
We define the $Kneser(k,l)$ graph, as usual. The graph is formed by taking $\binom{k}{l}$ vertices, each labeled by an $l$-subset of  $\big\{1,2, \ldots k \big\}$. There is an edge between two vertices if the two sets associated with them do {\bf not} share an element. For our purposes, we will usually be using the $Kneser(k,2)$ graph, which has $\binom{k}{2}$ vertices \cite[chap. 7]{Kn}. We shall denote the complete graph on $k$ vertices by $K_k$. For any vertex $j$ of $Kneser(k,2)$, we define $V_j$ to be the subset of $\big\{1,2,\ldots,k\big\}$ that is associated with $j$.
\\
Now, we define the $W$-Graph on $n=\binom{k+1}{2}$ vertices, $W(k)$, by the following adjacency matrix
$$
A(W(k))=\left[\begin{matrix} W_{11} & W_{12} \\ W_{21} & W_{22} \end{matrix} \right]
$$
where $W_{11}=A(K_k)$ , $W_{22}=A(Kneser(k,2))$ and $W_{12}$ is a $k-\text{by}-{k\choose 2}$ matrix such that:
$$
(W_{12})_{i,j} = \begin{cases}
1 & \text{if $i \in V_j$} \\
0 & \text{otherwise}
\end{cases}
$$
and let $W_{21} = W_{12}^{T}$. One can see that for $k \geq 4$ each row of $W_{12}$ has $k-1$ ones, which represent a maximal independent set of $Kneser(k,2)$. If $k \geq 5$ there are $k$ different maximal independent sets of size $k-1$, each represented by a different row.

This brings us to our next result:
\begin{theorem} \label{thm:W-graph} For all $k \ge 5$, $A(W(k))$ has exactly $k-1$ nonpositive eigenvalues.
\end{theorem}
\begin{proof}
It is known that $W_{11}$ has exactly $k-1$ nonpositive eigenvalues. We wish to show that $A(W(k))$ has the same number of nonpositive eigenvalues, so we examine the Schur complement $A(W(k))/W_{11} = W_{22} - W_{21}W_{11}^{-1}W_{12}$. Since by Lemma~\ref{lemm:schur}  $ i(W) = i(W_{11}) + i(A(W(k))/W_{11})$ , it is enough to show that $A(W(k))/W_{11}$ is positive definite.

By inspection, one can see that

$$
W_{11}^{-1} = \left[\begin{matrix}
-\frac{k-2}{k-1} & \frac{1}{k-1} &  \cdots & \frac{1}{k-1}\\
\frac{1}{k-1} & -\frac{k-2}{k-1} &  \cdots & \frac{1}{k-1}\\
\vdots & & \ddots & \vdots \\
\frac{1}{k-1} & \cdots & & -\frac{k-2}{k-1}\\
\end{matrix}\right]
$$

Note that each column of $W_{12}$ has exactly two entries which are equal to 1, and all the rest are zeros. Therefore,
$$
(W_{11}^{-1}W_{12})_{i,j} =\begin{cases}
-\frac{k-2}{k-1} + \frac{1}{k-1}=\frac{3-k}{k-1} & \text{if $i \in V_j$} \\
\frac{1}{k-1} + \frac{1}{k-1}=\frac{2}{k-1} & \text{otherwise.}\end{cases}
$$
By definition of $W_{21}$ we have
$$
(W_{21})_{i,j} = \begin{cases}
1 & \text{ if $j \in V_i$} \\
0 & \text{otherwise.}
\end{cases}
$$
Hence, the entries of $W_{21}W_{11}^{-1}W_{12} \in M_{{k\choose 2}}$ are as follows:

$$
(W_{21}W_{11}^{-1}W_{12})_{i,j} = \begin{cases}
\frac{3-k}{k-1}+\frac{3-k}{k-1}=\frac{6-2k}{k-1} & \text{if $V_i = V_j \Leftrightarrow i=j$}
\\
\frac{2}{k-1}+\frac{2}{k-1}=\frac{4}{k-1} & \text{if $V_i \cap V_j = \emptyset$}
\\
\frac{3-k}{k-1} + \frac{2}{k-1}=\frac{5-k}{k-1} & \text{otherwise.}
\end{cases}
$$

Comparing this to $W_{22}$ and using the definition $A(W(k))/W_{11} = W_{22}-W_{21}W_{11}^{-1}W_{12}$ we get that:

$$
\big(A(W(k))/W_{11}\big)_{i,j} = \begin{cases}
\frac{2k-6}{k-1}=1+\frac{k-5}{k-1} & \text{if $V_i = V_j \Leftrightarrow i=j$}
\\
1-\frac{4}{k-1}=\frac{k-5}{k-1} & \text{if $V_i \cap V_j = \emptyset$}
\\
\frac{k-5}{k-1} & \text{otherwise.}
\end{cases}
$$

Therefore, $A(W(k))/W_{11} = \frac{k-5}{k-1}J_{{k\choose 2}} + I_{{k\choose 2}}$, and its eigenvalues are $1$ with multiplicity ${k\choose 2}-1$ , and  $1+{k\choose 2}\frac{k-5}{k-1}$ with multiplicity $1$.
In conclusion, $A(W(k))/W_{11}$ is positive definite if and only if $k \ge 5$ and in this case $A(W(k))$ has exactly $k-1$ nonpositive eigenvalues.
\end{proof}
This has an immediate consequence for our bounds.
\begin{theorem} \label{thm:minbound} For all $k \ge 5$, $R(k,k+1) \ge NPO(k) > T_k$.
\end{theorem}
\begin{proof}
By Theorem~\ref{thm:W-graph}, there exists a graph, $W(k)$, on $n=\binom{k+1} {2}= T_k$ vertices whose adjacency matrix has only $k-1$ nonpositive eigenvalues. The upper bound follows from Theorem ~\ref{thm:rams}.
\end{proof}

\subsection{\textbf{Alternate Extreme Graphs}}
On 21 vertices, (we suspect that the exact value of $NPO(6)$ is 22), there are $4$ constructions similar to the $W$-graph. Instead of taking the complete graph on six vertices, we may take any two complete graphs that have six vertices together, so we get another three graphs on 21 vertices that has only 5 nonpositive eigenvalues.\\
Additionally, there is another graph of order 21, $G_6$, whose adjacency matrix has only $5$ nonpositive eigenvalues. We can construct this graph by taking two graphs - the \emph{Clebsch Graph}, which has 16 vertices, and $K_5$, the complete graph on 5 vertices, and connect each vertex from $K_5$ to a specific set of 8 vertices from the Clebsch Graph (this specific set forms an induced subgraph equal to $4K_2$). We get the graph $G_6$, a graph on 21 vertices that has only 5 nonpositive eigenvalues, which are\\
\big\{ -4.2361 ,-4.2361, -4.2361, -4.2361, -3.7346, 0.2361, 0.2361, 0.2361, 0.2361, 0.5853, 1, 1, 1, 1, 1, 1, 1, 1, 1, 1, 9.1493\big\}.

\begin{figure}[htb]
\centering
\caption{$G_6$, an extreme graph on $21$ vertices}
\includegraphics[width=0.55\textwidth]{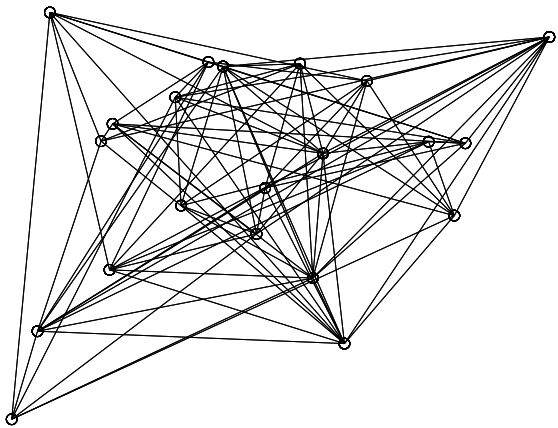}
\label{fig:twentyoneb}
\end{figure}
\pagebreak

\section{The Laplacian Matrix}
Having described the results found for adjacency matrices of graphs, we now turn to the Laplacian matrix. Most of the results found for the Laplacians are direct consequences of those found for the adjacency matrix. The Laplacian results were our original goal.
\begin{theorem} \label{thm:lap}
Let $k$ be a positive integer, and let $G$ a graph on $m$ vertices, $ m \ge NPO(k) $. Let $\{d_i\}$ be the degrees of $G$, arranged in non-increasing order and let $\{\lambda_i\}$ be the eigenvalues of the Laplacian matrix $L(G)$, arranged in non-increasing order. Then $\lambda_k \ge d_{NPO(k)}$.
\end{theorem}
\begin{proof}
We have $L(G)=D(G)-A(G)$, with $A(G)$ is the adjacency matrix of $G$ and $D(G)=diag\big(deg(v_1),deg(v_2)\ldots deg(v_m)\big)$. Let $H$ be a subgraph of $G$ induced by the $NPO(k)$ vertices with the largest degrees. Denote by $\hat{L}$ and $\hat{D}$ the submatrices of $L(G)$ and $D(G)$ respectively, that correspond to the vertices of $H$ (note that $\hat{L}$ is not the Laplacian matrix of $H$). By the Interlacing Theorem, $\lambda_k$ is bounded from bellow by the $k^{th}$ largest eigenvalue of $\hat{L}$. Since $H$ is of order $NPO(k)$, the $k^{th}$ largest eigenvalue of $-A(H)$ is nonnegative. Note that $\hat{L}=\hat{D}-A(H)$, and the result follows directly from Lemma~\ref{lemm:possemidef}.
\end{proof}
Thus we have the following results as consequences of Theorem~\ref{thm:lap} and earlier results on the adjacency matrix. In these corollaries, $\lambda_k$ shall refer to the $k^\text{th}$-largest eigenvalue of the Laplacian matrix $L(G)$.
\begin{corollary}
For any graph $G$ of order $m \ge 3$, $\lambda_2 \ge d_3$.
\end{corollary}
\begin{corollary}
For any graph $G$ of order $m \ge 6$, $\lambda_3 \ge d_6$.
\end{corollary}
\begin{corollary}
For any graph $G$ of order $m \ge 10$, $\lambda_4 \ge d_{10}$.
\end{corollary}
\begin{corollary}
For any graph $G$ of order $m \ge 16$, $\lambda_5 \ge d_{16}$.
\end{corollary}
In addition, we have a special corollary for regular graphs:
\begin{corollary}
Let $G$ be a $d$-regular graph on $n$ vertices, and let $k$ be a positive integer such that $NPO(k) \leq n < NPO(k+1)$. Then $L(G)$ has at least $k$ eigenvalues which are equal to or greater than $d$.\\
\end{corollary}

Theorem ~\ref{thm:lap} gives us the bound $\lambda_k \ge d_{NPO(k)}$. A natural question is, whether it is possible to improve it, i.e, is there exists a positive integer $m$ such that $m<NPO(k)$ and $\lambda_k \ge d_m$.
Note that Theorem~\ref{thm:lap} does not imply that such does not exist. For the case $k\leq 5$ we have examples of Laplacian matrices that show that the bounds given are the best possible bounds.\\

For $k \leq 3$ it follows directly from looking at the Laplacian matrices of the graphs $K_2$ and $C_5$. For $k=4$, take the graph $G_3$ (which appears in Figure ~\ref{fig:nineac}) and add four pendent vertices, each one of them is connected to one of the vertices of degree 3 in $G_3$ (there are four such vertices). This graph has 13 vertices, 9 of them of degree 4, but the fourth largest eigenvalue is smaller than 4, and hence this example shows that there exists a graph for which $\lambda_4 < d_9$. For $k = 5$, we start with $W(5)$. We add $30$ pendent vertices, and connect each one of them by an edge to $W(5)$ such that each vertex in the $Kneser(5,2)$ subgraph of $W(5)$ is connected by an edge to $3$ of these pendent vertices. The resulting Laplacian matrix has $\lambda_5 \approx 7.8438$ and $d_{15}=8$.
A similar process for larger $W$-graphs fails to work at some point, so that it is unclear whether the relationship on bounds given by Theorem~\ref{thm:lap} continues to be the best bound for Laplacian matrices beyond that point.

\section{Open Questions and Conjectures}
We defined $NPO(k)$ and determined the exact value of it for $k=1,2,3,4,5$. In addition, we gave an upper and lower bound on $NPO(k)$ for all positive integers $k$. The major question is: What is the exact value of $NPO(k)$ for each $k$?
\\
Moreover, the question of the bounds on Laplacian eigenvalues that motivated us remains open beyond $k=5$. Does there continue to be a precise relationship between the adjacency and Laplacian matrix? We know that the bound on the adjacency corresponds to a bound on the Laplacian. But does the \emph{best} bound on the adjacency correspond to the \emph{best} bound on the Laplacian?
\\
Finally, there may be other sorts of relationships, beyond just $k$ nonpositive eigenvalues. For example, for what sizes and what values of $l$ and $m$ can we say that the eigenvalues of the adjacency matrix must satisfy $\lambda_l+\lambda_m \le 0$? Such statements also translate to Laplacian eigenvalues bounds in terms of diagonal entries.
\\


\end{document}